\documentclass{amsart}
\usepackage{latexsym}
\usepackage[all]{xy}

\usepackage{amssymb} 
\usepackage{amsmath} 
\usepackage{amsthm}

\def\Z{{\mathbb{Z}}}
\def\K{{\mathbb{K}}}
\def\CC{{\mathbb{C}}}
\def\R{{\mathbb{R}}}
\def\C{{\mathbb{C}}}
\def\A{{\mathcal{A}}}
\def\B{{\mathcal{B}}}
\def\C{{\mathcal{C}}}
\def\F{{\mathbb{F}}}

\DeclareMathOperator{\codim}{codim}

\DeclareMathOperator{\Der}{Der}

\numberwithin{equation}{section}

\newcommand{\owari}{\hfill$\square$}

\newtheorem{theorem}{Theorem}[section]
\newtheorem{prop}[theorem]{Proposition}
\newtheorem{cor}[theorem]{Corollary}
\newtheorem{lemma}[theorem]{Lemma}
\newtheorem{define}[theorem]{Definition}
\newtheorem{rem}[theorem]{Remark}

\newtheorem{example}[theorem]{Example}

\title[Roots of characteristic polynomials]
{Roots of characteristic polynomials and 
intersection points of line arrangements}
\author{
Takuro Abe 
}

\address{Department of Mechanical Engineering and Science, Kyoto University, 
Kyoto 606--8501, Japan}
\email{abe.takuro.4c@kyoto-u.ac.jp}



\begin{document}

\maketitle

\begin{abstract}
We study a relation between roots of characteristic polynomials and 
intersection points of line arrangements. 
Using these results, we obtain a lot of applications for line arrangements. 
Namely, we give (i) a generalized addition theorem for line arrangements, 
(ii) a generalization of 
Faenzi-Vall\`{e}s' theorem  
over a field of arbitrary characteristic, (iii) 
a partial result on the conjecture of Terao for line arrangements, and (iv) 
a new sufficient condition for freeness over finite fields. Also, 
a higher dimensional version of our main results are considered. 
\end{abstract}

\section{Main results}
We use the notation in section two to state the main results in this article. 
Here some basic and special notations will be explained, which will be 
defined again in the next section. 

Let $\K$ be a field of arbitrary characteristic and consider affine line arrangements in $V=\K^2$. 
We say an affine line arrangement $\A$ is \textbf{free} with exponents 
$\exp_0(\A)=(d_1,d_2)$ if the cone $c\A$ of $\A$ is free with exponents 
$(1,d_1,d_2)$. 
For a line $H$, define 
$\A \cap H:=\{ H \cap H' \neq \emptyset \mid H' \in \A,\ H' \neq H\}$. 
Namely, this is the set of intersection 
points on $H$. Put $n_H:=|\A \cap H|$ and 
let $\chi(\A,t)$ be the characteristic polynomial of $\A$. 
Now let us state the main result in this article.

\begin{theorem}
Let $\C$ be an affine line arrangement
and assume that 
$\chi(\C,t)=(t-a)(t-a-b)$ with 
$a,b \in \mathbb{C}$ and $|a| \le |a+b|$. Then 
\newline
(1)\,\,
there are no $H \in \C$ such that $|a|<|\C \cap H|<|a+b|$. In other words, 
$\chi(\C,n_H) \ge 0$.
\newline
(2)\,\,
There are no line $L \not \in \C$ such that $|a|<|\C \cap L|<|a+b|$. In other words, 
$\chi(\C,n_L) \ge 0$.
\newline
(3)\,\,
Assume that $a,b \in \Z_{\ge 0}$. Then 
$\C$ is free if there is a line $H$ such that 
$|\C \cap H|=a$ or $a+b$. Equivalently, 
$\C$ is free if $\chi(\C,n_H)=0$ for some line $H$.
\label{elem}
\end{theorem}

If we assume the freeness, then we can obtain a stronger geometric condition on 
the arrangement.

\begin{cor}
In the same notation as in Theorem \ref{elem}, assume that 
$\C$ is free. Then 
\newline
(1)\,\,
$|\C \cap H| \in \Z_{\le a} \cup \{a+b\}$ for any $H \in \C$, and 
\newline
(2)\,\,
$|\C \cap L| \in \{a\} \cup \Z_{\ge a+b}$ for any line $L \not \in \C$.
\label{elem3}
\end{cor}

\begin{rem}
(1)\,\,
Theorem \ref{elem} (1) and (2) are non-trivial statements only when 
$a,b \in \R$ and $a<a+b$. 
\newline
(2)\,\,
Theorem \ref{elem} (1) gives some restriction on $H \in \C$ in terms of roots of $\chi(\C,t)$. 
On the other hand, Theorem \ref{elem} (2) seems to be more interesting. 
That is, the roots give a restriction on lines which are not belonging to $\C$. 
Hence Theorem \ref{elem} (2) says that combinatorics of $\C$ knows some information 
on geometry of $\C$. 
\newline
(3)\,\,
The case $n_H=a+b$ of Theorem \ref{elem} (3) when $b>0$ is essentially known to 
experts. See \cite{WY} for example. 
\newline
(4)\,\,
When $\K$ is a field of characteristic zero, 
Theorem \ref{elem} (3) and Corollary \ref{elem3} 
have been already proved by Faenzi and Vall\`{e}s as Propositions 5.2 and 
5.3 in \cite{FV}. 
Our proof is valid over a field of positive characteristic and 
the proofs in \cite{FV} and this article are very different. The former is algebro-geometric and 
in terms of Chern classes, and ours is algebraic and combinatorial, and in terms of Betti numbers.
\end{rem}

Let us check the statement in Theorem \ref{elem} and 
Corollary \ref{elem3} in the following example.

\begin{example}
(1)\,\,
The simplest but important example is a set of $n$-lines $\A$ in the 
real plane 
which go through the origin. Then it is obvious that 
$n_H=1$ for $H \in \A$, $n_L \in \{1,n-1,n\}$ for 
a line $L \not \in \A$ and $\chi(\A,t)=(t-1)(t-(n-1))$. 
This is trivial by using the property of parallel and generic lines, but 
Theorem \ref{elem} says that this holds true for all line arrangements.
\newline
(2)\,\,
Let $\A$ be an affine line arrangement in $\R^2$ defined by 
$$
x(x^2-y^2)(x^2-4y^2)(2x^2-y^2)(y-1)=0.
$$
Hence $|\A|=8$ and $$\chi(\A,t)=
t^2-8t+13=(t-
(4-\sqrt{3}))(t-
(4+\sqrt{3})).
$$
Hence Theorem \ref{elem} (1) and (2) say that 
$|H \cap \A| \neq 3,4,5$. In fact, we can check that 
$|H \cap \A| \in \{2, 7\}$ for $H \in \A$ and 
$|\A \cap L| \in \{1,2,6, 7,8\}$ for 
$L \not \in \A$.
\newline
(3)\,\,
Let $\A$ be an affine line arrangement in $\R^2$ defined by 
$$
xy(x^2-1)(y^2-1)(x^2-y^2)(x+y+1)(x+y-1)(x-y+1)(x-y-1)=0.
$$
Then $\chi(\A,t)=(t-5)(t-7)$, and it is easy to check that 
$|\A \cap H| =3$ or $5$ for any $H \in \A$, which matches Theorem \ref{elem} (1). 
Since we can check that there are no line $L 
\not \in \A$ such that 
$|L \cap \A|=6$, Theorem \ref{elem} (2) is satisfied. 
Also, Theorem \ref{elem} (3) 
shows that $\A$ is free.
\label{ex1}
\end{example}

The proofs of Theorem \ref{elem} and Corollary \ref{elem3} are simple, 
but we need algebraic methods for the proof of Theorem \ref{elem}. 
In particular, recent developments on exponents of two-dimensional multiarrangements 
(e.g., \cite{Y}, \cite{WY} and \cite{AN}) play 
the key roles. 

Recall that the coefficients of $\chi(\C,t)$ are  
the Betti numbers of the open manifold $V \setminus \cup_{H \in \C}H$ 
when $\K=\mathbb{C}$. 
Also, 
$\chi(\C,t)$ can be computed combinatorially in the arrangement cases. Hence 
we are interested in topological and combinatorial proofs of them.  
As far as we investigated, 
there are 
no such results similar to them.

Also, these results have a lot of applications. 
The first corollary is the following generalization of 
the addition theorem for line arrangements which includes 
a pair version of the conjecture by Terao. To state it, let us introduce some 
terminologies.
Define a \textbf{deletion pair of affine line arrangements} $(\A,\A')$ by 
$\A \supset \A'$ and $|\A'|+1=|\A|$. We say that a deletion pair $(\A,\A')$ is 
\textbf{free} if both $\A$ and $\A'$ are free. 
Then the following addition-type theorem holds.

\begin{cor}
A deletion pair $(\A,\A')$ is free if and only if 
$\chi(\A,t)$ and $\chi(\A',t)$ have a common root. In particular, the freeness of 
the deletion pair depends only on the combinatorics.
\label{ptc}
\end{cor}

Also, by using Theorem \ref{elem}, we can generalize 
Faenzi-Vall\`{e}s' theorem (Theorem \ref{main0}) in \cite{FV}. In Theorem \ref{main0}, 
the key condition is the existence of a point with multiplicity $h\ (n \le h \le n+r+1)$ 
for the arrangement $\A$ with $\chi(\A,t)=(t-n)(t-n-r)$. 
In this generalization, the role of this point is replaced by 
a free arrangement with exponents $(n-1,n-s)\ (s \ge 1)$, i.e., 
the following holds.

\begin{theorem}
Let $\K$ be a field of arbitrary characteristic and 
$\A$ an affine line arrangement such that $|\A|
=2n+r\ (n,r \in \Z_{\ge 0})$ and $\chi(\A,t)=(t-n)(t-n-r)$. Assume the following two 
conditions:
\begin{itemize}
\item[(1)]
$\A$ contains a free arrangement 
$\B$ with $\exp_0(\B)=(n-s,n-1)\ (s \ge 1)$, and 
\item[(2)]
there are no subarrangements $\C \subset \A$ such that $\B \subset \C$ and that 
$\chi(\C,t)=(t-n-u+1)(t-n+s)$ with $u > r+1$. 
\end{itemize}
Then 
$\A$ is free.
In particular, 
the freeness of such $\A$ depends only on combinatorics.
\label{main}
\end{theorem} 

We will explain in \S 4 why Theorem \ref{main} generalizes Faenzi-Vall\`{e}s' theorem. 

If we remove the assumption that ``$\B$ is free'' from the statement in Theorem \ref{main}, 
then can we say something on freeness and combinatorics? 
In fact, we can also show the following combinatorial 
statement on freeness.

\begin{theorem}
Let $\K$ be a field of arbitrary characteristic and 
$\A$ an affine line arrangement such that $|\A|
=2n+r\ (n,r \in \Z_{\ge 0})$ and $\chi(\A,t)=(t-n)(t-n-r)$. Assume that $\A$ contains an 
arrangement 
$\B$ with $\chi(\B,t)=(t-\alpha)(t-\beta)$ such that 
$\alpha \le \beta$ are real numbers with $\alpha \le n$ and $n-1 \le \beta$. Then 
$\A$ is free if and only if there is a line $H \in \A$ such that 
$n_H \in \{n,n+r\}$.
In particular, 
the freeness of such $\A$ depends only on combinatorics.
\label{mainc}
\end{theorem}


%

Also, we apply Theorem \ref{elem} and Corollary \ref{elem3} to 
obtain some results related to the conjecture of Terao (\S 5, 
Corollaries \ref{tc} and \ref{tc6}) and 
free arrangements over finite fields (\S 6, Theorem \ref{finite3}). 
Moreover, a higher dimensional version of Theorem \ref{elem} (1) and 
(2) will be given in \S 7.


The organization of this article is as follows. In \S 2
we introduce several definitions and results for the proof. Also, several lemmas for the proof of 
main results are proved. 
In \S 3 we prove main theorems.
In \S 4 we show generalized Faenzi-Vall\`{e}s' theorem as Theorems \ref{main} and \ref{mainc}. 
In \S 5 we show an application to the conjecture of Terao 
when one of the roots of the characteristic polynomial is at most five. In 
\S 6 we consider the case when the base field is a finite field. In \S 7, 
we give a higher dimensional version of Theorem \ref{elem}. 
\medskip




\noindent
\textbf{Acknowledgements}. 
This work is supported by JSPS Grants-in-Aid for Young Scientists (B) 
No. 24740012.

\section{Preliminaries}
In this section let us introduce several definitions and results, some of which have 
already defined in section one. We will use them throughout this article. 
We use \cite{OT} as a general reference on arrangement theory. Also, 
a recent paper \cite{Y3} is a nice reference on the algebraic aspects of multiarrangements.

Let $\K$ be a field of arbitrary characteristic unless otherwise 
specified, $V=\K^\ell$ and $S'=\mbox{Sym}^*(V^*) \simeq 
\K[x_1,x_2,\ldots,x_\ell]$ 
the coordinate ring of $V$. An \textbf{affine arrangement $\C$ of hyperplanes} in $V$ is 
a finite collection of affine hyperplanes in $V$.
Let $L(\C):=\{\cap_{H \in \B} H \neq \emptyset \mid \B \subset \C\}$ be the 
\textbf{intersection lattice} of $\C$. For $X \in L(\C)$, 
the \textbf{localization $\C_X$ of $\C$ at $X$} is defined by 
$$
\C_X:=\{H \in \C \mid X \subset H\}.
$$
Define 
$\mu:L(\C) \rightarrow \Z$ by $\mu(V)=1$, and by 
$\mu(X):=-\sum_{X \subsetneq Y \subset V} \mu(Y)$. Then the \textbf{characteristic polynomial} 
$\chi(\C,t)$ of $\C$ is defined by 
$$\chi(\C,t):=\sum_{X \in L(\C)} \mu(X) t^{\dim X}.
$$
For a hyperplane $H$, 
define $H \cap \C:=\{ H \cap H' \neq \emptyset \mid H' \in \C,\ H' \neq H\}$ and put 
$n_H:=|H \cap \C|$. Note that this definition is valid both when $H \in \C$ or $H \not \in \C$. 

Let $z$ be a new coordinate and 
define the \textbf{cone} $c\C$ of $C$ as follows. If $\C$ is defined by a non-homogeneous 
polynomial equation $Q=0$, then $c\C$ is defined by $z(cQ)=0$, where $cQ$ is the homogenized 
polynomial of $Q$ by the coordinate $z$. Hence $c\C$ is a \textbf{central} arrangement in $\K^{\ell+1}$, 
i.e., all hyperplanes contain the origin. We say that 
$\C$ is \textbf{essential} if $L(\C)$ contains the origin. 
For $H \in \C$, let $cH \in c\C$ denote the homogenized linear 
plane of $H$. 
Let 
$S:=\K[x_1,\ldots,x_\ell,z]$ and $\Der S$ be the module of $S$-derivations with 
a basis $\partial_{x_1},\ldots,\partial_{x_\ell},\partial_z$ dual to $x_1,\ldots,x_\ell,z$ respectively. 
Let $\alpha_{cH}$ be a defining linear 
form of $cH \in c\C$. Hence 
the defining polynomial $Q(c\C)$ of the cone $c\C$ of $\C$ is 
$z (\prod_{H \in \C} \alpha_{cH})$. 
Then define 
\begin{eqnarray*}
D(c\C):&=&\{\theta \in \Der S \mid \theta(\alpha_{cH}) \in S \alpha_{cH}\ (\forall H \in \C),\ 
\theta(z) \in S z\},\\
D_0(c\C):&=&\{\theta \in D(c\C) \mid \theta(z)=0\}.
\end{eqnarray*}
We say that $c\C$ is \textbf{free} with \textbf{exponents} $\exp(c\C)=(1,d_1,d_2,\ldots,d_\ell)$ 
if $D(c\C)$ is a free 
$S$-module with homogeneous basis elements $\theta_E=\sum_{i=1}^\ell x_i \partial_{x_i}+z\partial_z,\ 
\theta_1,\ldots, \theta_\ell$ such that $\deg \theta_i=d_i\ (i=1,\ldots,\ell)$. 
We say that an affine arrangement 
$\C$ is free with exponents $\exp_0(\C)=(d_1,\ldots,d_\ell)$ if $c\C$ is free with $\exp(c\C)=
(1,d_1,\ldots,d_\ell)$. 

Let $\A$ be a central hyperplane arrangement in $V$ 
and $m:\A \rightarrow \Z_{> 0}$ be a \textbf{multiplicity}. 
Here $\alpha_H$ denotes a defining linear form of $H \in \A$. 
Then a pair $(\A,m)$ is called a \textbf{multiarrangement} and we can define the 
logarithmic module 
$$
D(\A,m):=\{ \theta \in \Der S' \mid \theta(\alpha_H) \in S'\alpha_H^{m(H)}\ (\forall H \in \A)\}.
$$
Define 
$$
|m|:=\sum_{H \in \A} m(H).
$$
Let $Q(\A,m):=\prod_{H \in \A} 
\alpha_H^{m(H)}$. 
When $\ell=2$, $S'$ is two-dimensional. Hence 
$D(\A,m)$ is always free. Thus we can define its 
exponents $\exp(\A,m)=(d_1,d_2)$.

\begin{define}
For a central arrangement of hyperplanes $\C$ and $H_0 \in c\C$, 
let $(\C'',m)$ be the \textbf{Ziegler restriction} 
of $\C$ onto $H_0$ defined by $C'':=
\{H \cap H_0 \mid H \in \C \setminus \{H_0\}\}$
and by 
$$
m(H\cap H_0):=|\{ H' \in \C \setminus \{H_0\}  \mid H'\cap H_0=H \cap H_0\}|.
$$
The \textbf{Ziegler restriction of an affine arrangement $\C$ onto $H \in \C$} 
is that of $c\C$ onto $cH \in c\C$. Also, 
\textbf{Ziegler restriction of an affine arrangement $\C$ at infinity} 
is that of $c\C$ onto $\{z=0\} \in c\C$.
\label{z}
\end{define}

Now let us introduce a useful criterion for freeness. 

\begin{theorem}[Saito's criterion, \cite{Sa}, \cite{Z}]
Let $\theta_1,\ldots,\theta_\ell \in D(\A,m)$ 
be derivations with $\deg \theta_i=d_i\ (i=1,\ldots,\ell)$. Then they 
form a basis for $D(\A,m)$ if and only if $\theta_1,\ldots,\theta_\ell$ are $S'$-independent and 
$d_1+\cdots+d_\ell=|m|:=
\sum_{H \in \A}m(H)$.
\label{saito}
\end{theorem}

From now on, let us concentrate our interest on two-dimensional cases, i.e., 
line arrangements and its cones. 
Let $V=\K^2$ and $S'=\mbox{Sym}^*(V^*) \simeq 
\K[x,y]$ the coordinate ring of $V$. In this case, an affine  arrangement $\C$ in $V$ is 
a finite collection of affine lines in $V$. 
In this article, for a line arrangement $\C$, characteristic polynomial is denoted as follows:
$$
\chi(\C,t)=\sum_{X \in L(\C)} \mu(X) t^{\dim X}=t^2-|\C|t+b_2(\C).
$$
Here recall that $b_2(\C)$ coincides with the second Betti number of the open manifold 
$V \setminus \cup_{H \in \C} H$ when $\K=\CC$. 

Let $(d_1(\C),d_2(\C))$ denote the exponents of the Ziegler 
restriction of an affine line arrangement $\C$ onto $H \in \C$.
In general, we assume that 
$d_1(\C) \le d_2(\C)$. 
Then the following is the key theorem in this article.

\begin{theorem}[\cite{Y}, Theorem 3.2]
It holds that $\chi(\C,0)=b_2(\C) \ge d_1(\C)d_2(\C)$, and the equality holds if and only if $\C$ is free.
\label{yoshinaga}
\end{theorem}

Also, we use the results in the following papers; \cite{T}, \cite{T1}, 
\cite{Z}, \cite{WY}, \cite{AN} and \cite{A}. 
For the proof and application of main results, let us introduce some of them. 

First, let us introduce three results.  
Namely, the first one is the addition theorem in \cite{T}, the second 
the factorization theorem in \cite{T1}, and the third the Ziegler's restriction theorem 
in \cite{Z}. Note that all of these three were proved for any dimensional arrangements 
in these papers. Since we focus on line arrangements, we 
introduce 
the line arrangement cases of these results as follows.

\begin{theorem}[Addition theorem, \cite{T}]
Let $\A$ be an affine line arrangement and fix $H \in \A$. Define $\A':=\A \setminus \{H\}$ and 
$n_H:=|\A \cap H|$. Assume that $\chi(\A,n_H)=\chi(\A',n_H)=0$. Then $\A$ is free if and only if 
$\A'$ is free.
\label{addition}
\end{theorem}

\begin{theorem}[Factorization theorem, \cite{T1}]
Let $\A$ be a free affine line arrangement with $\exp_0(\A)=(d_1,d_2)$. Then 
$\chi(\A,t)=(t-d_1)(t-d_2)$.
\label{factorization}
\end{theorem}

\begin{theorem}[\cite{Z}]
If $\A$ is a free affine line arrangement with $\exp_0(\A)=(a,b)$, then 
its Ziegler restriction $(\A'',m)$ is free with $\exp(\A'',m)=(a,b)$. 
\label{ziegler}
\end{theorem}

The statements of Corollary \ref{ptc} and Theorem \ref{addition} 
are similar, and it is easy to see that 
the former is a generalization of the latter. 
The next two results are specialized ones for line arrangements. The first one 
is originally in \cite{WY}. 

\begin{lemma}[\cite{AN}, Lemma 4.2, Lemma 4.3]
Let $\A$ be a central line arrangement and $m,m'$ be multiplicities on $\A$ such that 
$|m|=|m'|+1$ and $m(H) \ge m'(H)$ for any $H \in \A$. If 
$\exp(\A,m')=(d_1,d_2)$, then $\exp(\A,m)=(d_1+1,d_2)$ or $(d_1,d_2+1)$.
\label{an1}
\end{lemma}

\begin{theorem}[\cite{A}]
Let $\A$ be an affine line arrangement defined over a field of 
characteristic zero. Put $\chi(\A,t)=(t-\alpha)(t-\beta)$ for $\alpha,\beta \in \mathbb{C}$. 
For the Ziegler restriction $(\A'',m)$ of 
$\A$ onto $H_0 \in \A$, put $\exp(\A'',m)=(d_1,d_2)$ with $d_1 \le d_2$. 
Assume that $ |m| \ge 2m(H) $ for any $H \in \A''$ and 
$|\A''|=:h>2$. Then 
\newline
(1)\,\,
$d_2-d_1 \le h-2$, and 
\newline
(2)\,\,
$||\alpha|-|\beta|| \le h-2$. In particular, 
$\A$ is free if $||\alpha|-|\beta|| \in \{h-2,h-3\}$. 
\label{h-2}
\end{theorem}

\begin{proof}
The statement (1) is the same as Theorem 3.5 in \cite{A}. 
Also, the statement (2) is essentially proved in 
\cite{A}. That is, combine $\Z \ni \alpha \beta \ge d_1 d_2$ (by Theorem \ref{yoshinaga}) and 
$\alpha+\beta=d_1+d_2=|\A|=|m|$ 
with (1) and Theorem \ref{yoshinaga}. 
\end{proof}
\medskip

When $(\A'',m)$ satisfies the condition $|m| \ge 2m(H)$ for any $H \in \A''$ 
in Theorem \ref{h-2}, we say that 
$(\A'',m)$ is \textbf{balanced}. 
We say that an affine line arrangement $\A$ is \textbf{balanced} if 
every Ziegler restriction of $\A$ is balanced. 
The following is famous in the theory of two-dimensional multiarrangements. 
We give a proof for the completeness.

\begin{lemma}
Let $\A$ be an affine line arrangement which is not 
balanced. 
Then the freeness of 
$\A$ depends only on $L(\A)$.
\label{balanced}
\end{lemma}

\begin{proof}
By definition, one of the Ziegler restrictions $(\A'',m)$ of $\A$ 
is not balanced. We may assume that $H:=\{x=0\} \in \A$ 
satisfies 
$2m(H) >|m|$. Let 
$\varphi:=(Q(\A'',m)/x^{m(H)}) \partial_{y}$. Then clearly 
$\varphi \in D(\A'',m)$ is a non-zero derivation 
of the smallest degree. 
Hence $\exp(\A'',m)$ is combinatorially determined as 
$(|m|-m(H),m(H))$ and Theorem \ref{saito} completes the proof. 
\end{proof}
\medskip

Now let us prove several statements for the proof of main results introduced in \S 1.
Some of them are well-known, but we give the whole proof for the completeness.  

\begin{lemma}
Let $\A$ be a central line arrangement with $|\A|=n$, 
$m$ be a multiplicity on $\A$ and put $\exp(\A,m)=(d_1,d_2)$ with 
$d_1 \le d_2$. 
\newline
(1)\,\,
If 
$|m| \ge 2n-2$, then 
$d_i \ge n-1$.
\newline
(2)\,\,
If 
$|m|\le 2n-2$, then 
$d_1 =|m|-n+1,\ d_2=n-1$.
\newline
(3)\,\,
Let $|m|=\alpha+\beta$ with $\alpha,\beta \in \R$ and $\alpha < \beta$. 
If $\alpha <n-1 < \beta$, then $\alpha < d_1 \le d_2 < \beta$.
\label{elem2}
\end{lemma}

\begin{proof}
(1)\,\,
Note that $\exp(\A)=(1,n-1)$. 
Take any multiplicity $m'$ such that $m(H)\ge  m'(H) \ge 1$ for any 
$H \in \A$ and $|m'|=2n-2$. 
Let $\theta_E$ be the Euler derivation. Then it is easily checked that 
$
\theta:=(Q(\A,m')/Q(\A)) \theta_E \in D(\A,m')$ is a 
non-zero element in $D(\A',m)$ of degree $n-1$ such that 
there are no $\theta' \in \Der S'$ satisfying $f \theta'=\theta$ for 
$f \in S'$ with $\deg f >0$.
%
Hence Theorem \ref{saito} implies that $\exp(\A,m')=(n-1,n-1)$. Since 
$D(\A,m') \supset D(\A,m)$, we complete the proof. 
\newline
(2)\,\,
Use the same $\theta=(Q(\A,m)/Q(\A)) \theta_E$ as in the proof of (1). 
Then $\deg \theta=|m|-n+1$ and 
it is clear that $\theta$ is a non-zero element of $D(\A,m)$ of 
the smallest degree. Hence 
Theorem \ref{saito} completes the proof. 
\newline
(3)\,\,
First assume that $\alpha \ge d_1$. Then 
the construction of $\theta$ in the proofs above shows that $d_2=n-1$. Hence 
$|m|=d_1+d_2 \le \alpha+n-1 < \alpha+\beta=|m|$, which is a contradiction. 
Hence 
$d_1 >\alpha$. Assume that $d_2 \ge \beta$. Then  
$|m|=d_1+d_2 > \alpha+\beta =|m|$, which is a contradiction. Hence $d_2 < \beta$. 
\end{proof}
\medskip

\begin{lemma}
Let $\A$ be a central line arrangement and $m,m'$ be multiplicities on $\A$ such that 
$m(H) \ge m'(H)$ for any $H \in \A$. Put $\exp(\A,m')=(d_1,d_2)$ and $\exp(\A,m)=(e_1,e_2)$ 
with $d_1 \le d_2,\ e_1 \le e_2$. Then $d_1 \le e_1,\ d_2 \le e_2$.
\label{an2}
\end{lemma}

\begin{proof}
Let $\theta_1,\theta_2\ ($resp:$\varphi_1,\varphi_2)$ be a basis for $D(\A,m')\ 
($resp:$D(\A,m))$ with $\deg \theta_i=d_i\ ($resp:$\deg \varphi_i=e_i)$. 
Since $D(\A,m) \subset D(\A,m')$, it is clear that $e_1 \ge d_1$. Assume that 
$e_2< d_2$. Then $\varphi_2=f \theta_1$ for $f \in S'$. 
Put $\varphi_1=g\theta_1+h\theta_2$ for $g,h \in S'$. Then 
the inequality $e_1 \le e_2<d_2$ shows that $h=0$. Hence $\varphi_1$ and 
$\varphi_2$ are $S'$-dependent, which is a contradiction. 
\end{proof}
\medskip


\begin{prop}
Let $\A \supset \B$ be affine line arrangements such that 
$\chi(\A,t)=(t-a)(t-c),\ \chi(\B,t)=(t-a)(t-b)$ with $a ,b,c \in \Z_{\ge 0}$. Assume that 
$a \le b \le c$. Then 
$\A$ is free if $\B$ is free.  
\label{elem5}
\end{prop}

\begin{proof}
Assume that $\B$ is free. Then $(d_1(\B),d_2(\B))=(a,b)$ by 
Theorem \ref{ziegler}.
By Theorem \ref{yoshinaga}, it suffices to show that $(d_1(\A),d_2(\A)) = (a,c)$. 
If not, then
Lemma \ref{an2} and Theorem \ref{yoshinaga} show a contradiction. 
\end{proof}
\medskip

\begin{example}
The inequality and the conditions on freeness in Proposition \ref{elem5} are essential. Consider 
\begin{eqnarray*}
\A:&=&xy(y^2-1)(x^2-4y^2)(x^2-9y^2),\\
\B:&=&x(y-1)(x^2-4y^2)(x^2-9y^2),\\
\C:&=&(x^2-4y^2)(x^2-9y^2).
\end{eqnarray*}
Then $\exp_0(\A)=(3,5),\ 
\exp_0(\C)=(1,3)$ and 
$\chi(\B,t)=(t-3)^2$, but $\B$ is not free.
\end{example}

\section{Proof of Theorem \ref{elem} and Corollary \ref{elem3}}

In this section we prove main results introduced in section one. 
\medskip

\noindent
\textit{Proof of Theorem \ref{elem}}. 
If both $a$ and $b$ are not real numbers, then 
$|a|=|a+b|$. Hence there is nothing to prove. 
So in the proof below, we may assume that $a$ and $b$ are both 
real numbers. Also, 
we may assume that $a$ and $a+b$ are both non-negative since 
the roots of $\chi(\C,t)=t^2-|\C|t+b_2(\C)$ are apparently non-negative.
Hence in the below, we may replace $|a|$ and $|a+b|$ by 
$a$ and $a+b$ respectively. 

(1)\,\,
Assume that such $H \in \C$ exists. Let $(\C'',m)$ be the Ziegler restriction of $\C$ onto 
$H$. Then $|C''|=n_H+1$ by definition of the cone. Thus 
$\exp(\C'')=(1,n_H)$ with $a<n_H<a+b$. Let $\exp(\C'',m)=(d_1,d_2)$ with $d_1 \le d_2$. 
Then it follows that $a<  d_1 \le d_2 < a+b$ by Lemma \ref{elem2} (3). 
Hence $d_1 d_2 >a(a+b)=b_2(\C)$, 
which 
contradicts Theorem \ref{yoshinaga}.
\medskip

(2)\,\,
First, note that the statement in (1) can be also written as 
$$
n_H^2-|\C|n_H+b_2(\C) \ge 0.
$$
Or equivalently, 
$$
b_2(\C) \ge n_H(|\C|-n_H).
$$
Now let $\B:=\C \cup \{L\}$. Then 
$|\B|=|\C|+1$ and $b_2(\C)=b_2(\B)+n_L$ by definition. Apply the inequality version 
of Theorem \ref{elem} (1) just above to $\B$ and $L \in \B$ to obtain 
$$
b_2(\B)=b_2(\C)+n_L \ge n_L(|\B|-n_L)=n_L(|\C|+1-n_L).
$$
Hence we conclude that 
$$
b_2(\C) \ge n_L(|\C|-n_L),
$$
which completes the proof. 

(3)\,\,
If $H \not \in \C$, replace $\C$ by $\C \cup \{H\}$ and we may assume that 
$H \in \C$ by Theorem \ref{addition}. 
First assume that $|\C \cap H|=a+b$. 
Let $(\C'',m)$ be the Ziegler restriction of $\C$ onto $H$. 
Then 
$\exp(\C'',m)$ is combinatorially determined as $(a,a+b)$ by Lemma \ref{elem2} (2). 
Hence $\C$ is free by Theorem \ref{yoshinaga}. 
Next assume that $|\C \cap H|=a$. Then $\exp(\C'')=(1,a)$. Hence Lemma \ref{elem2} (1) shows that 
$d_i \ge a$ for $\exp(\C'',m)=(d_1,d_2)$.
Again by Theorem \ref{yoshinaga}, 
we know that $a(a+b) \ge d_1d_2$. So Lemma \ref{an2} implies 
that $d_1=a,\ d_2=a+b$, which implies the freeness of 
$\C$ by Theorem \ref{yoshinaga}.\owari
\medskip

By Theorem \ref{elem} we know that $\chi(\C,n_H) \in \Z_{\ge 0}$. Since 
$\chi(\C,n_H)$ is a combinatorial invariant for $H \in \C$, it is natural to ask the 
meaning of this non-negative integer. The following is one of answers.

\begin{prop}
In the notation above, it holds that 
$$
\dim \mbox{coker}(\pi:D_0(c\C) \rightarrow D((c\C)'',m)) \le \chi(\C,n_H),
$$
where $((c\C)'',m)$ is the Ziegler restriction of $c\C$ onto $H \in \C$ and 
$\pi$ is the Ziegler restriction map (\cite{Z}).
\label{comb}
\end{prop}

\begin{proof}
Immediate from Theorem \ref{yoshinaga} and the proof of Theorem \ref{elem}.
\end{proof}
\medskip

\noindent
\textit{Proof of Corollary \ref{elem3}}.	
(1)\,\,
By Theorem \ref{elem}, it suffices to show that 
$|\C \cap H| \le a+b$ for $H \in \C$. Assume not. 
Then Lemma \ref{elem2} (2) shows that $d_1 d_2 <a(a+b)$ for 
$\exp(\C'',m)=(d_1,d_2)$, which contradicts Theorem \ref{yoshinaga}.  

(2)\,\,
First assume that $a=0$. This occurs only when 
all lines in $\C$ are parallel. In this case, Corollary \ref{elem3} is obvious. 
Hence we may assume that $a >0$. 

Since there is at least one point in $L(\C)$ by the 
previous paragraph, 
it holds that $\chi(\C,0)>0$ and $|\C| \ge 2$. Also, 
it is well-known 
that $\chi(\C,1) \ge 0$ (e.g., by Zaslavsky's theorem, \cite{Za}).  
Since $1 \le |\C|/2$, the non-negativity of $\chi(\C,0)$ and 
$\chi(\C,1)$ implies that $a \ge 1$.
Hence in the arguments below, we assume that $a \ge 1$. 

By Theorem \ref{elem}, it suffices to show that $|\C \cap L| \ge a$ for any line $L \not \in \C$. 
Assume not and put $C_1:=\C \cup\{L\}$. Let 
$(C_1'',m_1)$ be the Ziegler restriction of $\C_1$ at infinity and 
$n:=|C_1 \cap L| <a$. Then $b_2(\C_1)=b_2(\C)+n$. On the other hand, 
$\exp(\C_1'',m_1)=(a+1,a+b)$ or $(a,a+b+1)$ because  
$\exp(\C'',m)=(d_1,d_2)=\exp_0(\C)=(a,a+b)$ and Lemma \ref{an1}, where 
$(\C'',m)$ is the Ziegler restriction of $\C$ at infinity. 
Hence $a \ge 1$ implies that 
$b_2(\C_1)=b_2(C)+n =a(a+b)+n<a(a+b+1) \le (a+1)(a+b)$, which contradicts Theorem \ref{yoshinaga}. \owari
\medskip

\noindent
\textit{Proof of Corollary \ref{ptc}}. 
Let $a \in \mathbb{C}$ be a common root of $\chi(\A,t)$ and $\chi(\A',t)$. Recall the famous 
deletion-restriction 
formula 
$$
\chi(\A,t)=\chi(\A',t)-\chi(\A \cap H,t),
$$
where $\{H\}= \A \setminus \A'$. 
See \cite{OT}, Corollary 2.57 for example. Hence
$$
0=\chi(\A,a)=\chi(\A',a)-\chi(\A \cap H,a)=-\chi(\A \cap H,a), 
$$
By definition, $\chi(\A \cap H,t)=t-n_H$. Thus 
$\chi(\A \cap H,a)=a-n_H$. Hence $a=n_H \in \Z_{\ge 0}$, and 
both characteristic polynomials factorize into 
\begin{eqnarray*}
\chi(\A',t)&=&(t-a)(t-b),\\
\chi(\A,t) &=&(t-a)(t-b-1).
\end{eqnarray*}
Thus Theorem \ref{elem} (3) shows the freeness of both arrangements.\owari
\medskip

\begin{rem}
Corollary \ref{ptc} makes several proofs of the freeness of line arrangements 
easier, especially 
those related to extended Catalan and Shi arrangements. For example in \cite{A0}, 
the freeness of several deformations of the Coxeter arrangements of the type $A_2$ 
are proved by checking all the addition steps. However, if we use Corollary \ref{ptc}, 
it suffices to find a line $H$ on each deformations such that $n_H$ is one of the 
roots of their characteristic polynomials. 
\end{rem}

\begin{example}
Theorem \ref{elem} (3) and Corollary \ref{ptc} are useful as we saw above, but they are not 
enough to show freeness of all 
arrangements combinatorially. Recall the affine line arrangement $\A$ consisting of 
all edges and diagonals of a regular pentagon. Then 
$\chi(\A,t)=(t-5)^2$ and $\A$ is free, but $|\A \cap H|=4$ for any $H \in \A$. Hence 
we cannot apply Theorem \ref{elem} (3) and Corollary \ref{ptc} to show its freeness combinatorially. 
Of course, it is easy to see that there is a line $L \not \in \A$ such that $|\A \cap L|=5$. Hence 
Theorem \ref{elem} (3) shows that $\A$ is free, but this proof is not combinatorial. 
Also, it is easy to check that $\A$ contains a free arrangement with exponents $(3,3)$, but 
$\A$ does not satisfy the sufficient condition of freeness in Theorem \ref{main}. Hence 
the condition in Theorem \ref{main} is essential. 
\label{pentagon}
\end{example}

\begin{rem}
In the recent paper \cite{CH} by Cuntz and Hoge, an affine line arrangement which 
is free but is not obtained as the deletion pair is found.
\end{rem}

\section{Proof of Theorems \ref{main} and \ref{mainc}}
Before the proof of Theorem \ref{main} as an application of Theorem \ref{elem}, 
let us recall the following Faenzi-Vall\`{e}s' theorem.

\begin{theorem}[\cite{FV}, Theorem 2]
Let $\K$ be a field of characteristic zero and 
$\A$ be an affine $2$-arrangement in $V=\mathbb{K}^2$ such that $|\A|=2n+r\ (n,r \in \Z_{\ge 0})$ and 
that its localization $\B \subset \A$ at the origin consists of 
$h$-lines with $n \le h \le n+r+1$. If $\chi(\A,t)=(t-n)(t-n-r)$, then $\A$ is free.
\label{main0}
\end{theorem}


Now note that the following easy lemma.

\begin{lemma}
Let $\A$ be an affine line arrangement. Assume that $\A$ is free with 
$\exp_0(\A)=(n,n+r)\ (r \ge 0)$. Then there are no subarrangements $\B \subset \A$ with 
$\chi(\B,t)=(t-\alpha)(t-\beta),\ \alpha,\beta \in \Z$ and $\beta >n+r$.
\label{ez}
\end{lemma}

\begin{proof}
Assume that such $\B$ exists. By Theorem \ref{ziegler}, the Ziegler restriction of $\A$ 
at infinity is free with exponents $(n,n+r)$. Let $(d_1,d_2)$ be the exponents of the Ziegler restriction of 
$\B$ at infinity. Then Theorem \ref{yoshinaga} shows that $\alpha \beta  \ge d_1 d_2$. In other words,
$d_2  \ge \beta \ge n+r$, which 
contradicts Lemma \ref{an2}. 
\end{proof}
\medskip

If $\A$ contains a point with multiplicity $h$ with 
$n \le h \le n+r+1$, then it implies that 
$\A$ contains a free arrangement $\B$ with $\exp_0(\B)=(1,h-1)$ and 
$n-1 \le h-1 \le n+r$. Hence by applying Lemma \ref{ez}, we may regard 
Theorem \ref{main} as a generalization of Theorem \ref{main0} in the sense of freeness. 
Also, note that Theorem \ref{main} holds true over 
any fields of any characteristic. 

For the proof of Theorem \ref{main}, let us introduce the following 
corollary and lemma by using the results in the previous section. The 
first corollary might be similar to 
non-freeness criterion in \cite{K}. 

\begin{cor}
Let $\A \supset \B$ be an affine line arrangement such that 
$\chi(\A,t)=(t-a)(t-b),\ \chi(\B,t)=(t-c)(t-d)$ with integers  
$a \le b,\ c \le d$ and $\B$ is free. If $|\A \cap H| <b$, then 
for $H \in \A \setminus \B$, it holds that 
$
|\B \cap H|\in \{c\} \cup \{d,d+1,\ldots,a\}$.
\label{coro}
\end{cor}

\begin{proof}
Obvious by Theorem \ref{elem} and Corollary \ref{elem3}. 
\end{proof}
\medskip


\begin{lemma}
Let $\A$ and $\B$ be affine line arrangements such that 
$\A \supset \B$ with $|\A \setminus \B|=f$. Then we can order 
lines of $\A \setminus \B=\{H_1,\ldots,H_{f}\}$ in such a way that, 
for $\B_0:=\B,\ \B_i:=\B_{i-1} \cup \{H_i\}$ and 
$n_i:=|\B_{i-1} \cap H_i|$, it holds that 
$n_1\le n_2 \le \cdots \le n_f$.
\label{increase}
\end{lemma}

\begin{proof}
We use induction on $i$. First, 
let $H_1 \in \A \setminus \B$ be a line such that 
$|\B \cap H_1|=\min_{H \in \A \setminus \B}|\B \cap H|$. 
Then for any $H \in \A \setminus (\B \cup \{H_1\})$, it is obvious that 
$|\B \cap H_1| \le |(\B \cup \{H_1\}) \cap H|$. Assume that 
$H_1,\ldots,H_i \in \A \setminus \B$ satisfy the condition in the statement. 
Then choose $H_{i+1} \in \A \setminus \B_i$ such that 
$|\B_i \cap H_{i+1}|=\min_{H \in \A \setminus \B_i}|\B_i \cap H|$. 
Then it is obvious that $n_i \le |\B_i \cap H|$ for any 
$H \in \A \setminus \B_i$.
\end{proof}
\medskip

\noindent
\textit{Proof of Theorem \ref{main}}. 
If there is a line $H \in \A$ such that 
$n_H=n$ or $n+r$, then Theorem \ref{elem} shows that $\A$ is free. Assume that 
$n_H \neq n,\ n+r$. 
Again by Theorem \ref{elem}, we may assume that 
$n_H < n$ or $n_H >n+r$. Also, by Corollaries \ref{elem3} and \ref{coro}, 
$|H \cap \B|\in \Z_{\ge n-1} \cup\{n-s\}$ for $H \in \A \setminus \B$. 
%

Let $\A \setminus \B=\{H_1,\ldots,H_{r+s+1}\},\ \B_0:=\B$ and  
$\B_i:=\B \cup \{H_1\} \cup \cdots \cup\{H_i\}$. By Lemma \ref{increase}, we may assume that 
$n_1\le n_2 \le \cdots \le n_{r+s+1}$ for $n_i:=|\B_{i-1} \cap H_{i}|\ 
(i=1,\ldots,r+s+1)$. By the previous paragraph, 
we know that $ \{n-s\} \cup \Z_{\ge n-1} \ni n_1 \le n_{r+s+1} \in \Z_{<n} \cup \Z_{>n+r}$. 
Note that $n_{r+s+1}=|\A \cap H_{r+s+1}|$.

\noindent
\textbf{Case 1}. 
Assume that $n_1=n-s$. 

\noindent
\textbf{Case 1-1}. Assume that $n_2 >n-s$. Then 
$\B_1$ is free with $\exp_0(\B_1)=(n,n-s)$. By Theorem \ref{elem}, 
$n_2 \ge n$. Since $n \le n_2 \le n_{r+s+1} \neq n$, we have 
$n_{r+s+1} >n+r$ by Theorem \ref{elem}. Hence
$$
b_2(\A) > n(n-s)+(r+s-1)n+n+r=n(n+r)+r \ge n(n+r)=b_2(\A),
$$
which is a contradiction.

\noindent
\textbf{Case 1-2}. Assume that $n_1=\cdots=n_u=n-s<n_{u+1}$ for 
some $u>1$. 
Then $\B_u$ is free with $\exp_0(\B_u)=(n+u-1,n-s)$. 
%
If $r \ge u-1$, then $n_i \ge n+u-1 >n-1$ 
for $i>u$ by Corollary \ref{elem3} 
and $n+u-1 \le n+r$.
Hence 
\begin{eqnarray*}
b_2(\A) & >& (n+u-1)(n-s)+(r+s+1-u)(n+u-1)\\
&=&
(n+u-1)(n+r+1-u) \ge n(n+r)=b_2(\A)
\end{eqnarray*}
because of $0 \le r+1-u \le r$ and $n_{r+s+1}>n+r$, which is a contradiction.

If $r<u-1$, then there exists $\B \subset \C \subset \A$ such that 
$\chi(\C,t)=(t-n-u+1)(t-s)$ and $r<u-1$, which contradicts 
the condition (2). 

\noindent
\textbf{Case 2}. 
So we may assume that $n_1 \ge n-1$. 
If $n_{r+s+1}=n-1$, then 
$$
b_2(\A)=(n-1)(n-s)+(r+s+1)(n-1)=(n-1)(n+r+1)<n(n+r)=b_2(\A),
$$
which is a contradiction. Hence $n_{r+s+1} \ge n$. By the 
assumption and Theorem \ref{elem}, it holds that $n_{r+s+1} >n+r$. 
%
Hence 
$$
b_2(\A)>(n-1)(n-s)+(r+s)(n-1)+n+r=n(n+r)=b_2(\A),
$$
which is a contradiction. \owari
\medskip

It is natural to ask whether the same statement as in Theorem \ref{main} 
holds true for $s \le 0$. The answers is affirmative as follows.

\begin{prop}
In the same notation and condition as in Theorem \ref{main}, we assume that 
$-r \le s \le 0$. Then $\A$ is free if and only if $n_H \in \{ n,n+r\}$ for some 
$H \in \A$.
\label{mainprop}
\end{prop}

\begin{proof}
The ``if'' part follows by Theorem \ref{elem} (3). Conversely, assume that 
$\A$ is free and $n_H \not \in \{n,n+r\}$. Then Theorem \ref{elem} (1) shows that 
$n_H <n$ or $n_H>n+r$. Since $\A$ is free, Theorem \ref{yoshinaga} and 
Lemma \ref{elem2} (2) imply that $n_H<n$. 
Let $\A \setminus \B=\{H_1,\ldots,H_{r+s+1}\}$. Put 
$B_i$ and $n_i$ in the same way as in Theorem \ref{main} by Lemma \ref{increase}. 
Then Theorem \ref{elem} and Corollary \ref{elem3} show that 
$n-1 \le n_1 \le n_{r+s+1} \le n-1$. However, 
$$
b_2(\A)=(n-1)(n-s)+(r+s+1)(n-1)=(n-1)(n+r+1)<n(n+r)=b_2(\A),
$$
which is a contradiction. 
\end{proof}
\medskip




Before the proof of Theorem \ref{mainc}, we need the following lemma.

\begin{lemma}
Let $\A \supset \B$ be the same arrangements as in Theorem \ref{mainc}.  
Let us order $\A \setminus \B=\{H_1,\ldots,H_f\}\ (f:=2n+r-\alpha-\beta)$ 
in such a way that $\B_0:=\B,\ 
\B_i:=\B_{i-1} \cup \{H_i\}$ and 
$n_1 \le \cdots \le n_f$ for $n_i:=|\B_{i-1} \cap H_i|$ by 
Lemma \ref{increase}.
Let $a$ be the smallest integer 
satisfying $\alpha \le a$. Assume that $n_f \le n-1,\ 
n-1 < \beta$ and put 
$\chi(\B_i,t)=(t-\alpha_i)(t-\beta_i)$ with $|\alpha_i| \le |\beta_i|\ (i=1,\ldots,f)$. 
Then $\alpha_i$ and $\beta_i$ are both 
real numbers, and $\alpha_{i+1} \le \alpha_i \le \alpha \le 
\beta \le \beta_i \le \beta_{i+1}$ for any $i$. In particular,  
$n_i \le a$ for $i=1,\ldots,f$.
\label{a}
\end{lemma}

\begin{proof}
Let us prove by induction on $i$. Since 
$\chi(\B,t)=(t-\alpha)(t-\beta)$, Theorem \ref{elem} (1) shows the case 
$i=0$. Assume that the statement is true when $i \le k$. Since 
$n-1 < \beta \le \beta_k$, it holds that $n_{k+1}\le \alpha_k$ by Theorem \ref{elem} (2). 
Since 
$$
\chi(\B_{k+1},t)=t^2-(\alpha_k+\beta_k+1)t+\alpha_k\beta_k+n_{k+1}, 
$$
the roots of this polynomial are of the form
$$
t_{\pm}=\displaystyle \frac{\alpha_k+\beta_k+1 \pm \sqrt{(\alpha_k-\beta_k-1)^2+4(\alpha_{k}-n_{k+1})}}{2}.
$$
Since $\alpha_k \ge n_{k+1}$, it follows that $t_{\pm} \in \R$. 
Also, it is easy to see that 
$
t_- \le \alpha_k$ and $\beta_k \le t_+$. 
Hence Theorem \ref{elem} (1) completes the proof. 
\end{proof}
\medskip

\noindent
\textit{Proof of Theorem \ref{mainc}}. 
The ``if'' part follows from Theorem \ref{elem} (3). Assume that 
$\A$ is free and there are no $H \in \A$ such that $n_H \in \{n,n+r\}$. 
By Lemma \ref{balanced} we may assume that $n_H \le n+r$. Hence 
Theorem \ref{elem} (1) shows that $n_H \le n-1$ for $H \in \A$. 
Let us order $\A \setminus \B=\{H_1,\ldots,H_f\}\ (f:=
2n+r-\alpha-\beta)$ in such a way that $\B_0:=\B,\ 
\B_i:=\B_{i-1} \cup \{H_i\}$ and 
$n_1 \le \cdots \le n_f \le n-1$ for $n_i:=|\B_{i-1} \cap H_i|$ by 
Lemma \ref{increase}. 



\textbf{Case 1}. 
Assume that $\alpha \le n-1 \le \beta$. If $\beta=n-1$, then 
\begin{eqnarray*}
b_2(\A)&\le&\alpha(n-1)+(n-1)(2n+r-n+1-\alpha)\\
&=&(n-1)(n+r+1)<n(n+r)=b_2(\A),
\end{eqnarray*}
which is a contradiction. Hence we may assume that $n-1 < \beta$. 

Let $a,b$ be integers such that $\alpha \le a < \alpha+1$ and 
$\beta-1<b \le \beta$. Hence $\alpha+\beta=a+b$. Since 
$\alpha +\beta =|\A| \in \Z$, it holds that $a \le n-1 \le  b$ and 
$\alpha \beta \le ab$. 
Since $n-1 < \beta$, we may apply 
Lemma \ref{a} to obtain that 
$n_i \le a$.  
Hence 
\begin{eqnarray*}
b_2(\A)&\le&ab+a(2n+r-a-b)\\
&=&a(2n+r-a)<n(n+r)=b_2(\A),
\end{eqnarray*}
which is a contradiction. 

\textbf{Case 2}. Assume that $n-1<\alpha \le \beta < n$. 
Then 
$\alpha+\beta =2n-1$ and 
$\alpha \beta \le (n-\displaystyle \frac{1}{2})^2$. 
Hence 
\begin{eqnarray*}
b_2(\A)&\le&(n-\displaystyle \frac{1}{2})^2+(n-1)(2n+r-2n+1)\\
&=&(n-1)(n+r+1)+\displaystyle \frac{1}{4}<n(n+r)=b_2(\A),
\end{eqnarray*}
which is a contradiction. 

\textbf{Case 3}. Assume that $n-1 < \alpha \le n,\  n \le \beta $. 
Let $a$ and $b$ be the same integers as in the Case 1. Hence $n \le b$ and 
$a =n$. 
Since $n_i \le n-1$ and $n \le \beta$, we may apply 
Lemma \ref{a} to obtain that 
$n_i \le a$. Hence 
\begin{eqnarray*}
b_2(\A)&\le&nb+n(2n+r-n-b)\\
&=&n(n+r)=b_2(\A).
\end{eqnarray*}
The equality holds only when $\alpha=n=n_1=\cdots=n_f$,  
which contradicts $n_f \le n-1$. \owari

%

\medskip




\section{Applications related to the conjecture of Terao}
In this section we study the relation between the conjecture of Terao and 
the results in the previous sections. 

First, let us show the following proposition, which is a generalization of 
Theorem \ref{main} in a special case.

\begin{prop}
Let $\A$ be an affine line arrangement such that $\chi(\A,t)=(t-n)(t-n-r)$ 
with $n \in \Z_{\ge 0}$ and $r \in \Z_{\ge 1}$. Assume that 
$\A$ contains an arrangement $\B$ with $\chi(\B,t)=(t-n+2)^2$. 
Then $\A$ is free if and only if $n_H=n$ or $n+r$ for some $ H \in \A$.
\label{n-2}
\end{prop}

\begin{proof}
The ``if'' part follows by Theorem \ref{elem} (3). Assume that 
$\A$ is free and $n_H \not \in \{n,n+r\}$. 
Then $n_H>n+r$ or $n_H <n$ by Theorem \ref{elem}
Also, $n_H >n+r$ implies 
the non-freeness of $\A$ by Theorem \ref{yoshinaga} and Lemma \ref{elem2} (2). 
Hence we may assume that $n_H <n$. 

Let $\{H_1,\ldots,H_{r+4}\}=\A \setminus \B$. Put 
$\B_0:=\B,\ \B_i:=\B_{i-1} \cup \{H_i\}\ (i=1,\ldots,r+4)$. Then 
for $n_i:=|H_i \cap \B_{i-1}|$, we may assume that $n_1\le n_2\le \cdots\le n_{r+4}<n$ by 
Lemma \ref{increase}. 
Then 
\begin{eqnarray*}
b_2(\A) &\le& (n-2)^2+(n-1)(2n+r-(2n-4))\\
&=&n(n+r)-r <n(n+r)=b_2(\A),
\end{eqnarray*}
which is a contradiction.
\end{proof}
\medskip
 
Theorem \ref{n-2} has the following corollary. 

\begin{cor}
Let $\A$ be an affine line arrangement.

\noindent
(1)\,\,
If $\chi(\A,t)=(t-2)(t-2-r)$ with $r>0$, then 
the freeness of $\A$ depends only on $L(\A)$. 
\newline
(2)\,\,
If $L(\A)$ contains a point and $\chi(\A,t)=(t-3)(t-3-r)$ with $r > 0$, then 
the freeness of $\A$ depends only on $L(\A)$.
\label{coro2}
\end{cor}

\begin{proof}
(1)\,\,
Since $\A$ contains an empty arrangement with exponents $(0,0)$, Proposition \ref{n-2} 
completes the proof.
\newline
(2)\,\,
Since $\A$ contains a Boolean arrangement with exponents $(1,1)$, Proposition \ref{n-2} completes 
the proof. 
\end{proof}
\label{coro3}
\medskip

The following can be proved by the same way as in Theorem \ref{n-2}.

\begin{prop}
Let $\A$ be an affine line arrangement such that $\chi(\A,t)=(t-n)(t-n-r)$ 
with $n,r\in \Z_{\ge 0}$.

\noindent
(1)\,\,
Assume that $r \ge 2$ and 
$\A$ contains an arrangement $\B$ with $\chi(\B,t)=(t-n+2)(t-n+3)$. 
Then $\A$ is free if and only if $n_H=n$ or $n+r$ for some $ H \in \A$.

\noindent
(2)\,\,
Assume that $r \ge 4$ and 
$\A$ contains an arrangement $\B$ with $\chi(\B,t)=(t-n+3)^2$. 
Then $\A$ is free if and only if $n_H=n$ or $n+r$ for some $ H \in \A$.
\label{n-3}
\end{prop}

On the conjecture of Terao, which asserts that the freeness of an arrangement $\A$ 
depends only on its combinatorics $L(\A)$, we can give a few contribution by using these with 
Theorem \ref{h-2}. The conjecture of Terao for line arrangements in $\mathbb{C}^2$ 
is confirmed when $|\A| \le 10$ by Wakefield-Yuzvinsky (\cite{WY}, Corollary 7.5), and $|\A| \le 11$ 
by Faenzi-Vall\`{e}s.  
(\cite{FV}, Theorem 5). 

Now using the results in this article, first, we can show the following.

\begin{cor}
Let $\A$ be an affine line arrangement in $\mathbb{C}^2$ such that 
$\chi(\A,t)=(t-n)(t-n-r)$ with $n,r\in \Z_{\ge 0}$. If 
$r \ge n-3$, then the freeness of $\A$ depends only on $L(\A)$.
\label{tc}
\end{cor}

\begin{proof}
Let $(\A'',m)$ be the Ziegler restriction of $\A$ at infinity. 
By Lemma \ref{balanced}, we may assume that 
$\A$ and $(\A'',m)$ are balanced. 
Put 
$\exp(\A'',m)=(d_1,d_2)$ with $d_1 \le d_2$. 
By Theorem \ref{h-2} (2), we know that the combinatorial invariant  
$h:=|\A''| \ge r+2$. 
When 
$h=r+2$ or $r+3$, the freeness of $\A$ is confirmed by Theorem \ref{h-2} (2). 
Assume that $h \ge r+4 \ge n+1$. Then Theorem \ref{elem} (1) shows that 
$h \not \in \{n+2,\ldots,n+r\}$, and Theorem \ref{elem} (3) shows 
that $\A$ is free when $h=n+1$ or $n+r+1$. 
Also, the non-freeness of $\A$ when $h >n+r+1$ is checked in \cite{WY}, or 
by applying Theorem \ref{yoshinaga} and Lemma \ref{elem2} (2). 
\end{proof}
\medskip

Using Corollary \ref{tc}, in this article, 
we check the conjecture of Terao from a different point of view from 
\cite{WY} and \cite{FV}. Namely, we prove the conjecture under the restriction on 
the roots of characteristic polynomials, not on  
the number of lines.

\begin{cor}
Let $\A$ be an affine line arrangement in $\mathbb{C}^2$ such that 
$\chi(\A,t)=(t-n)(t-n-r)$ with $n,r \in \Z_{\ge 0}$. If 
$\{ n, n+r\} \cap \{0,1,2,3,4,5\} \neq \emptyset$, then 
the freeness of $\A$ depends only on $L(\A)$.
\label{tc6}
\end{cor}

\begin{proof}
If $\{n,n+r\} \cap \{0,1\} \neq 0$, then the conjecture of Terao is easy to check. 
Assume that $n+r \in \{2,3,4,5\}$. Then \cite{WY} and \cite{FV} complete the proof. 
So we may assume that $n \in \{2,3,4,5\}$. Also, the case $r=0$ can be verified 
by \cite{WY} and \cite{FV}. So assume that $r>0$.

Assume that $n=2$. Then Corollary \ref{coro2} (1) completes the proof. Assume that 
$n=3$. Then a point is contained in $L(\A)$. Hence Corollary \ref{coro2} (2) completes the proof. 

Assume that $n=4$. By Lemma \ref{balanced}, 
we may assume that $\A$ is balanced. 
Then Corollary \ref{tc} verifies the statement 
when $r \ge 1$. Hence it suffices to check when $\chi(\A,t)=(t-4)^2$, 
which is checked in \cite{WY} and \cite{FV}.

Assume that $n=5$.
By Lemma \ref{balanced}, 
we may assume that $\A$ is balanced. 
Then Corollary \ref{tc} verifies the statement 
when $r \ge 2$. Hence it suffices to check when $\chi(\A,t)=(t-5)(t-6)$ or 
$(t-5)^2$, 
which is checked in \cite{FV}. 
\end{proof}

\section{The case over finite fields}
In this section let us consider the case when $\K$ is a finite field $\F_q$. 
We give an another proof of Theorem 10 in \cite{Y2}. Also, we give 
a new sufficient condition for freeness which is 
a similar result to that in \cite{Y2}. Namely, in \cite{Y2}, it is shown that
an arrangement which has $q$ as the root of the characteristic polynomial is free. 
Here we show that the same holds true when $q-1$ is a root.  

In this section we use the following 
setup. 
Let $\F_q$ be a finite field of cardinality $q=p^n$ for a prime number $p$ 
and $V =\F_q^2$. Recall that, for an affine line arrangement $\A$ in $V$, it holds that 
$$
\chi(\A,q)=|V \setminus \cup_{H \in \A}H|.
$$
See Theorem 2.69 in \cite{OT} for example. Now 
consider 
a multiarrangement $(\A,m)$ in $V$. Put $\exp(\A,m)=(d_1,d_2)$ with $d_1 \le d_2$.

\begin{prop}
Assume that $m(H) \le q$ for any $H \in \A$. Then 
\newline
(1)\,\,
the inequality $d_1 < q<d_2$ cannot occur.
\newline
(2)\,\,
If $|m| \ge 2q$, then $d_1=q$.
\newline
(3)\,\,
If $|m| = 2q-1$, then $d_2=q$.
\label{finite}
\end{prop}

\begin{proof}
(1)\,\,
Let $\theta_1,\theta_2$ be a basis for $D(\A,m)$ with $\deg \theta_i=d_i$. 
Assume that $d_1 < q <d_2$. 
Note that $\varphi:=x^q \partial_x +y^q \partial_y \in D(\A,m)$, which is of degree $q$. 
Hence $ \varphi=f \theta_1$ for some polynomial $f$. Since 
$\varphi$ has no divisors in $\Der (S')$, this is a contradiction.
\newline 
(2)\,\,
By (1) and $|\A|=d_1+d_2 \ge 2q$, we know that 
$d_2 \ge d_1 \ge q$. 
Since $\varphi \in D(\A,m)$, 
we know that $d_1 \le q$
, which completes the proof. 
\newline
(3)\,\,
By assumption, $d_2 \ge q$. If $d_2>q$, then $d_1 < q < d_2$, which is a contradiction.
\end{proof}
\medskip

The following is proved in \cite{Y2}. Here we give an another proof of it.

\begin{cor}[\cite{Y2}, Theorem 10]
Let $\A$ be an affine line arrangement in $V$. 
\newline
(1)\,\,
If $\chi(\A,q)=0$, then $\A$ is free.
\newline
(2)\,\,
If $|\A| \ge 2q-1$ and $\A$ is free, then $\chi(\A,q)=0$. 
\label{finite2}
\end{cor}

\begin{proof}
Let $(\A'',m)$ be the Ziegler restriction of 
$\A$ at infinity. Put $\exp(\A'',m)=(d_1,d_2)$ with $d_1 \le d_2$. Note that 
$d_1+d_2=|\A|$. Also, 
note that we may apply Proposition \ref{finite} since 
the base field is $\F_q$. 

(1)\,\,
Let $\chi(\A,t)=(t-q)(t-r)$. Note that $q+r=d_1+d_2 =|\A|=|m|$. 
First assume that $r \le q$. Then Theorem \ref{yoshinaga} implies that 
$qr \ge d_1d_2$. Hence $d_1 \le r \le q \le d_2$. By Proposition \ref{finite} (1), 
we know that $q=d_1$ or $q=d_2$. Hence $\A$ is free by Theorem \ref{yoshinaga}. 

Second assume that $r>q$. Then again the inequalities $d_1 \le q <
r \le  d_2$ and Proposition \ref{finite} (1) show that 
$d_1=q$, which implies the freeness. 

(2)\,\,
Since $|m|=|\A| \ge 2q-1$, Proposition \ref{finite} (2) and (3) imply
that $d_1=q$ or $d_2=q$. Then the freeness of $\A$, Theorems \ref{factorization} and 
\ref{ziegler} 
complete the proof. 
\end{proof}
\medskip

By applying Theorem \ref{elem}, we can prove 
the following new result on arrangements in $\F_q^2$. 

\begin{theorem}
Let $\A$ be an affine arrangement in $V=\F_q^2$.
If $\chi(\A,q-1)=0$, then $\A$ is free.
\label{finite3}
\end{theorem}

\begin{proof}
Put $\chi(\A,t)=(t-q+1)(t-q+r)$ with $r \in \Z$. 
Since $\chi(\A,q)=r=|V \setminus \cup_{H \in \A} H| \ge 0$, 
we know that $r \in \Z_{\ge 0}$, and 
$\A$ is free if $r=0$ by Corollary \ref{finite2}. Assume that 
$r \ge 1$. 
Since $\chi(\A,0) \ge 0$, it holds that 
$\chi(\A,q)=r \le q$. Let $V \setminus \cup_{H \in \A} H=\{
p_1,\ldots,p_r\}$ and we may assume that $p_1$ is the origin. 
Then there are $(q+1)$-lines containing $p_1$ and not belonging to $\A$. Hence 
there is at least one line $L \not \in \A$ such that $p_1 \in L$ and 
$p_i \not \in L$ for 
$i=2,\ldots,r$. 
Then $|\A \cap L|=q-1$. Hence Theorem \ref{elem} (3) shows that 
$\A$ is free. 
\end{proof}

\section{Higher dimensional version}

In this section we prove a higher dimensional version of Theorem \ref{elem}. 
Unless otherwise specified, we use the following notation in this section. 
Let $\A$ be an 
affine arrangement of hyperplanes in $V=\K^\ell$ with $\ell \ge 3$. Let 
$L_i(\A):=\{X \in L(\A) \mid \codim_VX=i\}$ and denote 
$\chi(\A,t)=t^\ell-b_1t^{\ell-1}+b_2t^{\ell-2}+\cdots+(-1)^\ell b_\ell$. When $\K=\CC$, 
$b_i$ is the $i$-th Betti number of the open manifold 
$V \setminus \cup_{H \in \A}H$. It is known that 
$b_1=|\A|$. 

The following is a direct generalization of Theorem \ref{elem} to an arbitrary dimensional 
arrangements.

\begin{theorem}
For a hyperplane $H$, let $\chi(\A \cap H,t)=\sum_{i=0}^{\ell-1} (-1)^{i}c_i t^{\ell-1-i}$. If we 
put $|\A \cap H|=:h$, then 
$$
b_2 \ge c_2+(b_1-h-1)h.
$$
In particular, when 
$(b_1-1)^2-4b_2+4c_2 \ge 0$, there are no hyperplanes $L$ such that 
$$
\displaystyle \frac{b_1-1-\sqrt{(b_1-1)^2-4b_2+4c_2}}{2}<
|\A \cap L|
<
\displaystyle \frac{b_1-1+\sqrt{(b_1-1)^2-4b_2+4c_2}}{2}.
$$
\label{ineq}
\end{theorem}

To prove Theorem \ref{ineq}, let us recall one definition and 
introduce two results. 

\begin{prop}[\cite{AY}, Theorem 4.1 (1)] 
Let $(\B,m)$ the Ziegler restriction of $c\A$ onto $H \in c\A$. 
Define $b_2(\B,m):=
\sum_{X \in L_2(\B)} d_1^Xd_2^X$, where 
\begin{eqnarray*}
\B_X:&=&\{H \in \B \mid X \subset H\},\\
m_X:&=&m|_{\B_X},\\
\exp(\B_X,m_X):&=&(d_1^X,d_2^X,0,\ldots,0).
\end{eqnarray*}
Then $b_2(\A) \ge b_2(\B,m)$.
\label{b2}
\end{prop}

\begin{proof}
Let us recall the definition of the  
characteristic polynomial 
$\chi(\B,m,t)$ of the multiarrangement $(\B,m)$ (Definition 2.6, \cite{ATW}).
Then the local-global formula (Theorem 3.3, \cite{ATW}) shows that 
the $b_2(\B,m)$ above coincides with that of the coefficient of $t^{\ell-2}$ of 
$\chi(\B,m,t)$. Hence the inequality is nothing but 
Theorem 4.1 (1) in \cite{AY}.
\end{proof}
\medskip

\noindent
\textit{Proof of Theorem \ref{ineq}}. 
Assume that $H \in \A$. Let $c\A$ be the cone of $\A$ and 
$(\B,m)$ be the Ziegler restriction of $c\A$ onto $cH$. 
By Proposition \ref{b2}, we know that $b_2(\A) \ge b_2(\B,m)$.

Next assume that $H \not \in \A$. Let $c\A_1$ be the cone of $\A_1:=\A \cup \{H\}$ and 
$(\B_1,m_1)$ be the Ziegler restriction of $c\A_1$ onto $cH$. 
By Proposition \ref{b2}, we know that $b_2(\A \cup \{H\})=b_2+h \ge b_2(\B_1,m_1)$. 

Hence it suffices to show that 
$b_2(\B,m) \ge c_2+h(b_1-h-1)$ when $H \in \A$, and 
$b_2(\B_1,m_1) \ge c_2+h(b_1+1-h-1)$ when $H \not \in \A$. Since $|\A \cap H|=h$, we know that $|\B|=
|\B_1|=h+1$. 
Hence the following Lemma \ref{p1} completes the proof. \owari
\medskip

\begin{lemma}
Let $(\B,m')$ be a multiarrangement and $m_H:\B \rightarrow \{0,1\}$ be the multiplicity defined by 
$m_H(L):=\delta_{H,L}$ for $H, L \in \A$. Define $m:=m'+m_H$ and let 
$|\B|-1=:h$. Then 
$$
b_2(\B,m) \ge b_2(\B,m')+h.
$$
\label{p1}
\end{lemma}

\begin{proof}
For $X\in L_2(\B)$, define 
$\exp(\B_X,m_X)=:(d_1^X,d_2^X,0,\ldots,0)$ and
$\exp(\B_X,m_X')=:(e_1^X,e_2^X,0,\ldots,0)$. Then 
Proposition \ref{b2} shows that 
$$
b_2(\B,m)-b_2(\B,m')=\sum_{X \in L_2(\B),\ X \subset H} 
(d_1^Xd_2^X-e_1^Xe_2^X).
$$
Recall that $\exp(\B_X)=(1,|\B_X|-1,0,\ldots,0)$. Hence 
Lemmas \ref{an1} and \ref{elem2} show that 
$$
d_1^Xd_2^X -e_1^X e_2^X \ge |\B_X|-1.
$$
Since $\sum_{H \supset X \in L_2(\B)} (|\B_X|-1)=h$, it holds that 
$$
b_2(\B,m)-b_2(\B,m') \ge h,
$$
which completes the proof. 
\end{proof}

\begin{rem}
Theorem \ref{ineq} can be also proved by applying Theorem \ref{elem} (1) and (2) 
with the combinatorial restriction map in \cite{AY}.
\end{rem}

\begin{example}
Let $\overline{\A}$ be a Weyl arrangement of the type $B_4$ defined by 
$$
xyzw(x^2-y^2)(x^2-z^2)(x^2-w^2)(y^2-z^2)(y^2-w^2)(z^2-w^2)=0
$$
and 
$\A:=\overline{\A}|_{\alpha_H=1}$ for some $H \in \A$. Then 
$\chi(\A,t)=(t-3)(t-5)(t-7)$. Also, 
$\chi(\A \cap H,t)=(t-1)(t-3)(t-5)$ for any $H \in \A$. Hence 
$b_1=15,\ b_2=71$ and $c_2=23$ in the notation of Theorem \ref{ineq}. Hence Theorem \ref{ineq} shows that 
there are no $L$ such that $|\A \cap L|=7$. 
\end{example}

Theorem \ref{ineq} is not easy to apply. To make it useful, let us prove the following Lemma. 

\begin{lemma}
Let $\B$ be an essential arrangement in $V$ with $|\B|=h+1$. Then 
$b_2(\B) \ge (\ell-1)(h-\ell+2)+(\ell-1)(\ell-2)/2$.
\label{el}
\end{lemma}

\begin{proof}
We use the double induction on $\ell$ and $h$. When $\ell=1$ there is nothing to show.
Note that the essential arrangement in $V=\K^\ell$ requires $|\B| \ge \ell$. 
When $h+1=\ell$, $\B$ is nothing but the Boolean arrangement. Hence it is free 
with exponents $(1,\ldots,1)$. Thus $b_2(\B)=\ell-1+(\ell-1)(\ell-2)/2$. 

Now let $\B$ be an arbitrary essential arrangement in $V=\K^\ell$ with $|\B| >\ell$. 
Then obviously there is a hyperplane $H \in \B$ such that $\B':=\B \setminus \{H\}$ and 
$\B'':=\B \cap H$ are both essential. Now apply the induction assumption to obtain 
that 
\begin{eqnarray*}
b_2(\B') &\ge&  (\ell-1)(h-\ell+1)+(\ell-1)(\ell-2)/2,\\
b_1(\B'') &\ge&  \ell-1.
\end{eqnarray*}
By the deletion-restriction formula which appeared in the proof of Corollary \ref{ptc}, 
we know that 
$$
b_2(\B)=b_2(\B')+b_1(\B''),
$$
which completes the proof. 
\end{proof}
\medskip

\begin{cor}
For a hyperplane $H$, let $h:=|\A \cap H|$. Assume that 
$d_{\A,H}:=(b_1+\ell-2)^2-4b_2(\A)-2(\ell-1)(\ell-2) \ge 0$ and 
$L_\ell(\A) \neq \emptyset$. Then 
there are no hyperplane $H$ such that 
$$
\displaystyle \frac{b_1+\ell-2 - \sqrt{d_{\A,H}}}{2} < 
| \A \cap H| <
\displaystyle \frac{b_1+\ell-2 + \sqrt{d_{\A,H}}}{2}.
$$
In particular, it holds that 
$$
h^2-b_1h+b_2 \ge 0.
$$
\label{ineq2}
\end{cor}

\begin{proof}
Combine Theorem \ref{ineq} and Lemma \ref{el}.
\end{proof}
\medskip

\begin{rem}
When $\ell=2$, Corollary \ref{ineq2} is nothing but Theorem \ref{elem} (1) and 
(2). 
\end{rem}

\begin{example}
Let $\A$ be an affine arrangement of planes in $V=\K^3$ defined by 
$$
(x\pm 1)(x\pm 2)(x \pm 3)(x \pm 4)y(y \pm 1)(z \pm 1)=0.
$$
Then it is easy to check that $\chi(\A,t)=(t-2)(t-3)(t-8)$. Hence 
$|\A|=13$ and $b_2(\A)=46$. Since $L_3(\A) \neq \emptyset$, 
Corollary \ref{ineq2} shows that there are no planes $L$ such that 
$$
6 \le |\A \cap L| \le 8.
$$
\end{example}


\begin{thebibliography}{OSS}

 \bibitem{A0} T. Abe, The stability of the family of $A_2$-type arrangements. 
 \textit{J. Math. Kyoto Univ} \textbf{46} (2006), no. 3, 617--636.

 \bibitem{A} T. Abe, Exponents of 2-multiarrangements and freeness of 3-arrangements. 
 \textit{J. Alg. Combin.} \textbf{38} (2013), no. 1, 65--78.

 \bibitem{AN} T. Abe and Y. Numata, 
 Exponents of $2$-multiarrangements and multiplicity lattices. 
\textit{J. Alg. Combin.} \textbf{35} (2012), no. 1, 1--17.


 \bibitem{ATW} T. Abe, H. Terao and M. Wakefield, 
The characteristic polynomial of a multiarrangement. 
\textit{Adv. in Math.} \textbf{215} (2007), 825--838. 

 \bibitem{AY}
 T. Abe and M. Yoshinaga, 
 Free arrangements and coefficients of characteristic polynomials. 
\textit{Math. Z.} \textbf{275} (2013), 911--919.

 \bibitem{CH}
 M. Cuntz and T. Hoge, 
 Free but not recursively free arrangements. To appear in \textit{Proc. Amer. Math. Soc.} (2013).  

 \bibitem{FV}
 D. Faenzi and J. Vall\`{e}s, Logarithmic bundles and Line arrangements, an approach via the standard construction. 
arXiv:1209.4934v1.

 \bibitem{K}
 J. Kung, A geometric condition for a hyperplane arrangement to be free. 
 \textit{Adv. in Math.} \textbf{135} (1998), 303--329. 


 \bibitem{OT} P. Orlik and H. Terao, \textit{Arrangements of hyperplanes}.
 Grundlehren der Mathematischen Wissenschaften, 
 \textbf{300}. Springer-Verlag, Berlin, 1992.

 \bibitem{Sa} {K. Saito},
 Theory of logarithmic differential forms and logarithmic vector fields. 
 \textit{J. Fac. Sci. Univ. Tokyo Sect. IA  Math}. 
 \textbf{27} (1980), 265--291. 



 \bibitem{T} H. Terao, Arrangements of hyperplanes and their freeness I, II. 
 \textit{J. Fac. Sci. Univ. Tokyo} \textbf{27} (1980), 293--320.   

 \bibitem{T1} H. Terao, 
 Generalized exponents of a free arrangement of hyperplanes and
 Shepherd-Todd-Brieskorn formula.  \textit{Invent. Math.}  {\bf 63}  
(1981), no. 1,
159--179. 

\bibitem{WY} M. Wakefield and S. Yuzvinsky, Derivations of an effective divisor 
on the complex projective line. \textit{Trans. Amer. Math. Soc.} \textbf{359} (2007), 
no. 9, 4389--4403. 


 \bibitem{Y} M. Yoshinaga, On the freeness of 3-arrangements. 
 \textit{Bull. London Math. Soc.} \textbf{37} (2005), no. 1, 126--134. 

 \bibitem{Y2} M. Yoshinaga, Free arrangements over finite fields. 
 \textit{Proc. Japan Acad. Ser. A} \textbf{82} (2006), no. 10, 179--182. 
 
 \bibitem{Y3} M. Yoshinaga, Freeness of hyperplane arrangements and related topics. 
 To appear in \textit{Annales de la Faculte des Sciences de Tolouse}. arXiv:1212.3523. 



\bibitem{Za}
T. Zaslavsky, 
Facing up to arrangements: Face-count formulas for partitions of spaces by hyperplanes.
\textit{Memoirs Amer. Math. Soc.}, \textbf{154}, 1975.  


\bibitem{Z}
G. M. Ziegler, 
Multiarrangements of hyperplanes and their freeness. in
{\it Singularities} (Iowa City,
IA, 1986), 345--359, Contemp. Math., {\bf 90}, Amer. Math. Soc.,
Providence, RI, 1989.

\end{thebibliography}
\end{document}